\documentclass[reqno, intlimits, sumlimits]{amsart} 

\usepackage{graphicx,amsmath, amstext, amssymb, amsthm, bbm, lmodern, enumerate, esint}
\usepackage{hyperref}
\hypersetup{colorlinks=true,linkcolor=black,citecolor=black,filecolor=black,urlcolor=black}

	\renewcommand{\Re}{\text{Re}}			
	\renewcommand{\Im}{\text{Im}}			
	
    \renewcommand{\tilde}{\widetilde}
    \renewcommand{\theta}{\vartheta} 
    \renewcommand{\phi}{\varphi}
    \renewcommand{\epsilon}{\varepsilon}
    
	\newcommand{\lebesgue}		{\ensuremath{\lambda\mkern-8mu\lambda}}
	\newcommand				{\eins}			{\mathbbm{1}}   
	\newcommand				{\norm}[1]		{\left\lVert#1\right\rVert}
	\newcommand				{\abs}[1]		{\left\lvert#1\right\rvert}

	\DeclareMathOperator	{\IC}			{\mathbb{C}}
	\DeclareMathOperator	{\IE}			{\mathbb{E}} 
	\DeclareMathOperator	{\IN}			{\mathbb{N}}
	\DeclareMathOperator	{\IP}			{\mathbb{P}}
	
	\DeclareMathOperator	{\IR}			{\mathbb{R}}
	
	\DeclareMathOperator	{\Var}			{\mathbb{V}\text{ar}}
	\DeclareMathOperator	{\supp}			{supp}

\theoremstyle{plain}
\newtheorem{theorem}			{Theorem}
\newtheorem{lemma}	[theorem]	{Lemma}
\newtheorem{corollary}	[theorem]	{Corollary}
\newtheorem{proposition}	[theorem]	{Proposition}
                     
\theoremstyle{definition}
\newtheorem{definition}	[theorem]	{Definition}

\begin{document}
 \title[Rate of Convergence to the Circular Law]{Rate of Convergence to the Circular Law via Smoothing Inequalities for Log-Potentials}
\author{Friedrich G\"otze}
\address{Friedrich G\"otze, Faculty of Mathematics, Bielefeld University, Germany}
\email{goetze@math.uni-bielefeld.de}
\author{Jonas Jalowy}
\address{Jonas Jalowy, Faculty of Mathematics, Bielefeld University, Germany}
\email{jjalowy@math.uni-bielefeld.de}
\date{\today}
\thanks{Supported by the German Research Foundation (DFG) through the IRTG 2235}
 \begin{abstract}
 The aim of this note is to investigate the Kolmogorov distance of the Circular Law to the empirical spectral distribution of non-Hermitian random matrices with independent entries. The optimal rate of convergence is determined by the Ginibre ensemble and is given by $n^{-1/2}$. A smoothing inequality for complex measures that quantitatively relates the uniform Kolmogorov-like distance to the concentration of logarithmic potentials is shown. Combining it with results from Local Circular Laws, we apply it to prove nearly optimal rate of convergence to the Circular Law in Kolmogorov distance. Furthermore we show that the same rate of convergence holds for the empirical measure of the roots of Weyl random polynomials.
 \end{abstract}
 
 \maketitle 
 
\section{Introduction}\label{sec:Intro}
The (complex) empirical spectral distribution of a non-Hermitian random matrix with i.i.d. entries will converge to the uniform distribution on the complex disc as the size of the matrix tends to infinity. This \emph{Circular Law} has a long history going back to Ginibre \cite{Gin65}, proving the special case of complex Gaussian entries. Later, Bai \cite{bai97} used Girko's Hermitization Trick, introduced in \cite{Girko85}, to prove the Circular Law under extra density and moment assumptions. The density assumption was removed by G\"otze and Tikhomirov \cite{GT07} and several reductions of the moment conditions appeared in \cite{GT10Circular,PZ10,TV08Circular}. Significant progress was possible due to the control of the smallest singular values in \cite{Rudelson08,TV09Discrete,RV08littlewood,TV08Circular}. Ultimately, the Circular Law was proven under optimal second moment assumption by Tao and Vu (with an appendix by Krishnapur) \cite{TV10Circular}. We recommend the survey \cite{BC12} for further discussions.

Like the history of the Circular Law already indicates, non-Hermitian Random Matrix Theory has become a fast growing field with recent activity in the last years. Applications are various and include dynamics of (neural) networks, scattering in chaotic quantum systems and Coulomb plasma, see \cite{AR06NeuralNetwork,FKS97,AB11,For10} just to name a few. Random Matrix Theory is mostly concerned with universality phenomena, like the global universality in the Circular Law. Here, the limiting spectral distribution remains universal among a big class of entry distributions of the underlying matrix. Its local analogue has recently been investigated in \cite{BYY14Local,BYY14LocalEdge,GNT17Local,TV15uni,AEK19} among others. In this note, we address universality of the rate of convergence, containing local as well as global universality in a uniform and quantitative manner.

Consider a non-Hermitian random matrix $X=(X_{ij})_{1\leq i,j\leq n}$ having independent real or complex entries $X_{ij}$, where in the complex case we additionally assume $\Re X_{ij}$ and $\Im X_{ij}$ to be independent. Define the empirical spectral distribution by 
\begin{align*}
\mu_n=\frac{1}{n}\sum_{j=1}^n \delta_{\lambda_j (X/\sqrt n)}, 
\end{align*}
where $\delta_\lambda$ are Dirac measures in the eigenvalues $\lambda_j$ of the scaled matrix $X/\sqrt n$.
The Circular Law states that if $\IE X_{ij}=0$ and $\IE \abs{X_{ij}}^2=1$, then $\IP$-a.s. we have
\begin{align*}
 \mu_n\Rightarrow\mu_\infty\text{, where }d\mu_\infty(z)=\frac 1\pi\eins_{B_1(0)}(z)dz
\end{align*}
is the uniform distribution on the complex disc. We abbreviate $dz=d\lebesgue(z)$ for the Lebesgue measure $\lebesgue$ on $\IC$.

We are interested in the rate of convergence, more precisely in the Kolmogorov distances over balls 
\begin{align*}
D_n:=D(\mu_n,\mu_\infty):=\sup_{z_0\in\IC,R>0}\abs{\mu_n(B_R(z_0))-\mu_\infty (B_R(z_0)) } 
\end{align*}
as $n\to\infty$. In \cite{GS02}, $D$ is called Discrepancy metric. The study of Kolmogorov-like metrics of complex measures is widely uncommon in the literature of non-Hermitian Random Matrix Theory so far. Therefore, let us provide some additional information about advantages of studying $D_n$. Most importantly, convergence in this distance coincides with weak convergence in the case of an absolutely continuous limit distribution, see Lemma \ref{lem:Convergence}. Hence $D$ is a reasonable object to study the rate of convergence to the Circular Law, in particular because some explicit calculations of $D_n$ are possible. Using the rotational symmetry of $\mu_\infty$ and the mean empirical spectral distribution $\bar\mu _n =\IE \mu_n$ of the so called Ginibre ensemble, i.e. $X_{ij}\sim\mathcal N_{\IC}(0,1)$, an elementary calculation shows
\begin{lemma}\label{lem:GinibreRate}
The mean empirical spectral distribution $\bar\mu _n =\IE \mu_n$ of the Ginibre ensemble satisfies
\begin{align}\label{eq:GinibreRate}
D(\bar\mu_n,\mu_\infty)\sim \frac 1{\sqrt{2\pi n }}
\end{align}
and
\begin{align}\label{eq:GinibreRateBulk}
\sup_{\substack{B_R(z_0)\subseteq \IC\setminus B_{1+\epsilon}(0)\\ \text{or }B_R(z_0)\subseteq B_{1-\epsilon}(0)}}\abs{\bar\mu_n(B_R(z_0))-\mu_\infty (B_R(z_0)) }\lesssim e^{-n\epsilon^2}.
\end{align}
\end{lemma}

Here and in the sequel $\sim$ denotes asymptotic equivalence, $\lesssim$ will denote an inequality that holds up to a parameter-independent constant $c>0$ that may differ in each occurrence. Moreover we write $A\asymp B$ if $c\abs B\le\abs A\le C\abs B$ for some constants $0<c<C$.

According to \eqref{eq:GinibreRate}, the optimal rate of $\mu_n$ to the Circular Law turns out to be $\mathcal O (1/\sqrt n)$. 
Interestingly, if one avoids the edge of $B_1(0)$ by a fixed distance $\epsilon$, then \eqref{eq:GinibreRateBulk} implies that the rate of convergence is exponentially fast.

Nevertheless we cannot expect an exponentially fast rate of convergence for the \emph{non-averaged} empirical spectral distribution $\mu_n$, since it is still sensitive to individual eigenvalue fluctuations. In particular, for each fixed set of eigenvalues $\{ \lambda_i\}_{i\le n}$ we may select a ball of radius $(10\sqrt n)^{-1}$ contained in $B_1(0)$ such that it does not cover any eigenvalue and obtain $D_n\gtrsim 1/n$. Heuristically, the typical distance of $n$ uniformly distributed eigenvalues is $n ^{-1/2}$. Therefore one may vary $B_R(z_0)$ up to a magnitude of $n^{-1/2}$ without covering a new eigenvalue and hence we expect $D_n$ to be of order $n ^{-1/2}$ as well. Our main result, Theorem \ref{thm:rateofconv} below, states that nearly the optimal rate of convergence still holds for non-Gaussian entry distributions of the underlying matrix. If all moments of the entry distribution are finite, we will show that
\begin{align}\label{eq:rateofconvwop}
D_n\lesssim n^{-1/2+\epsilon}
\end{align}
holds with overwhelming probability. A sequence of events $\Omega_n$ is said to hold with overwhelming probability (in short w.o.p.) if $\IP(\Omega_n^c)\lesssim n^{-Q}$ for any $Q>0$.

\begin{figure}[h]\label{fig:Samples}
 \includegraphics[width=.32\textwidth]{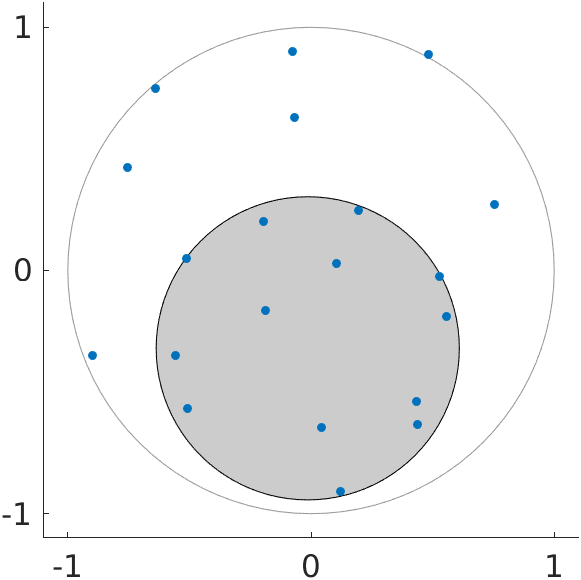}\hfill\includegraphics[width=.32\textwidth]{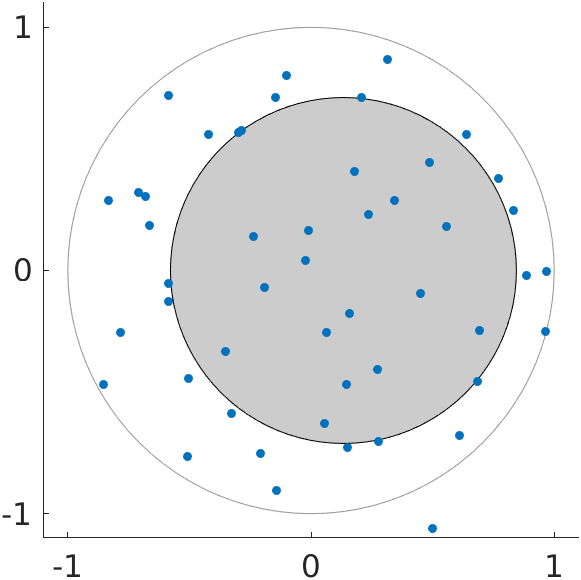}\hfill\includegraphics[width=.32\textwidth]{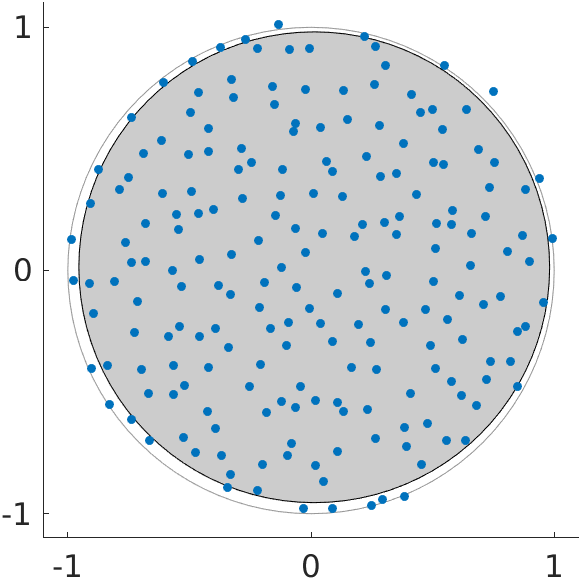}
  \caption{Samples of the spectrum of $X$ for $n=20$, $50$ and $200$ and a gray Ball $B_R(z_0)$ that attains the supremum in $D_n$. Theorem \ref{thm:rateofconv} below shows that clusters (like in the left sample) and sparse areas (like in the middle sample) do not significantly differ from the limit distribution. Here, we chose the entries to be uniformly distributed over a centered square in $\IC$. Even for these non-Gaussian entries, we clearly expect $B_R(z_0)$ to be close to $B_1(0)$ for larger $n$. This statement is exact for Ginibre matrices, as we see in the proof of Lemma \ref{lem:GinibreRate}.}
\end{figure}

In the remaining part of this section, we will introduce logarithmic potentials, so that the statement of the smoothing inequality in section \ref{subsec:Smoo} will be more clear. Different variants of results concerning the rate of convergence will be given in section \ref{subsec:roc}. After a discussion of the results, we will also discuss random polynomials in section \ref{sec:Polynomials}. Sections \ref{sec:Smoothing} and \ref{sec:Rates} contain the proofs of our results.

Similar to the role of the Stieltjes transform in the theory of Hermitian random matrices, the weak topology of measures $\mu$ on $\IC$ can be expressed in terms of the so called logarithmic potential $U$, which is the solution of the distributional Poisson equation. More precisely for every finite Radon measure $\mu$ on $\IC$ the \emph{logarithmic potential} defined by
\begin{align}\label{eq:logPot}
 U_\mu(z):=-\int_{\IC}\log\abs{t-z}d\mu (t)=(-\log\abs \cdot\ast\mu)(z)
 \end{align}
 satisfies the equation $\Delta U=-2\pi \mu $ in the sense of distributions. Obviously the logarithmic potential of a measure is superharmonic in $\IC$, harmonic outside the support of $\mu$ and is only unique up to addition of harmonic functions. The advantage of the logarithmic potentials $U_n$ of $\mu_n$ in non-Hermitian random matrix theory is the following identity known as \emph{Girko's Hermitization trick}
\begin{align}
U_n(z)&=-\frac{1}{n}\sum_{j=1}^n\log\abs{\lambda_j-z}=-\frac{1}{n}\log\abs{\det \Big(\frac 1 {\sqrt n}X -z\Big)} \nonumber\\
&=-\frac{1}{n}\log\det \sqrt{\Big(\frac 1 {\sqrt n}X -z\Big)\Big(\frac 1 {\sqrt n}X -z\Big)^*} =-\int_0^\infty\log (x)d\nu_n^z(x),\label{eq:GirkoHerm}
\end{align}
where $\nu_n^z$ is the empirical singular value distribution of the shifted matrix $X/\sqrt n -z$. Due to this fact, all the information on the complex spectrum of $X/\sqrt n$ is stored in the real and positive spectra of $(X/\sqrt n -z)(X/\sqrt n -z)^*$ for \emph{all} shifts $z$. Note that its symmetrized version around $0$ is the empirical eigenvalue distribution of the Hermitian matrix
\begin{align*}
 V(z)=\begin{bmatrix}
  0 & (X/\sqrt n -z)\\
  (X/\sqrt n -z)^* & 0
 \end{bmatrix}.
\end{align*}
 Under certain conditions on the matrix entries, the logarithmic potential $U_n$ concentrates around the logarithmic potential $U_\infty$ of the Circular Law given by
\begin{align*}
 U_\infty(z)=\begin{cases}
              -\log\abs z &\text{, if }\abs z >1\\
              \tfrac 1 2 (1-\abs z ^2) &\text{, if }\abs z \le 1
             \end{cases}.
\end{align*}
Let us fix some notation and the above-mentioned conditions.
\begin{definition}\label{def:cond}
\begin{enumerate}
 \item[(A)] We say $X$ satisfies \emph{condition (A)} if it has independent entries $X_{ij}$ with mean zero, variance $\IE \abs{X_{ij}}^2=1$ and if for all $p\in \IN$ it holds $\max_{i,j}\IE\abs{X_{i,j}}^p<\infty$.
\item[(B)] We say $X$ satisfies condition (B) if it has independent entries, where 
$$\max_{i,j}\abs{\IE X_{ij}}\le n^{-1-\epsilon}\text{ and }\max_{i,j}\abs{1-\IE \abs{X_{ij}}^2}\le n^{-1-\epsilon}$$ for some $\epsilon>0$ and furthermore 
 $$\max_{i,j,n} \IE \abs{X_{ij}}^{4+\delta}<\infty$$ for some $\delta>0$.\end{enumerate}
\end{definition}

Note that in contrast to Wigner matrices, the distributions of the entries may be different and clearly, (A) implies (B). 

%

Results on the concentration of the logarithmic potentials are used to derive Local Circular Laws. This has been explicitly proven in \cite[Theorem 25]{TV15uni} for subexponential entries, which match Gaussian moments up to third order. A very recent Local Circular Law by Alt, Erd\H{o}s and Kr\"uger \cite{AEK19} allows to extract the following concentration of logarithmic potentials from their proof, which weakens these assumptions.
 
\begin{proposition}[\cite{AEK19}]\label{prop:ConcLogPotAEK}
If $X$ obeys (A), then for every $\epsilon,\tau,Q>0$ there exists a constant $c>0$ such that
\begin{align}
\IP\left(\abs{U_n(z)-U_\infty(z)}\leq cn^{-1+\epsilon}\right)\geq 1-n^{-Q}
\end{align}
holds for any $z\in B_{1+\tau}(0)$.
\end{proposition}
Note that the result in \cite{AEK19} holds in a more general setting of inhomogeneous variances under some additional assumptions. Therefore, it is possible to show a rate of convergence in Kolmogorov distance $D$ for the inhomogeneous Circular Law as well. However, we stick to normalized variances in this work in order to avoid exhaustive notation and technicalities.

G\"{o}tze, Naumov and Tikhomirov \cite{GNT17Local} further weakened the assumptions, improved the rate and the result has been generalized to products of independent matrices, but at the cost of restricting the
region to the bulk $\abs{\abs z -1}\geq \tau$. 

\begin{proposition}[\cite{GNT17Local}]\label{thm:ConcLogPot}
If $X$ obeys (B), then for every $\tau,Q>0$ there exist a constant $c>0$ such that 
\begin{align}
\IP\left(\abs{U_n(z)-U_\infty(z)}\leq c\frac{\log^4 n}{n}\right)\geq 1-n^{-Q}
\end{align}
holds for any $z\in B_{1+\tau^{-1}}(0)$ such that $\abs{1-\abs z }\ge \tau\}$.
\end{proposition}

Since Proposition \ref{prop:ConcLogPotAEK} and Proposition \ref{thm:ConcLogPot} are not explicitly worked out in \cite{GNT17Local} and \cite{AEK19} respectively, we will derive them in Appendix \ref{sec:Appendix} based on the results proved in the corresponding paper.

\section{Main Results}\label{sec:Main}
\subsection{Smoothing Inequality}\label{subsec:Smoo} Consider a sequence of probability measures $\mu_n$ on $\IC$ with logarithmic potentials $U_n$. If $U_n$ converges pointwise to some function $U:\IC\to(-\infty,\infty]$ and if $U_n$ is locally uniformly Lebesgue integrable, then (by continuity of $\Delta$ on the space of distributions) there exist a probability measure $\mu=-\frac{1}{2\pi}\Delta U $ on $\IC$ such that $\mu_n$ converges weakly to $\mu$. The following smoothing inequality quantifies this statement by relating $D_n$ to the concentration of logarithmic potentials.

\begin{theorem}\label{prop:Smoothing}
Let $\mu,\nu$ be probability measures on $\IC$ with $\supp\nu\subseteq B_K(0)$ for some $K>0$, let $U_\mu, U_\nu$ be their logarithmic potentials and fix some $1\le p\leq \infty$. For any $a\ge 1/2$ we have
\begin{align*}
D(\mu,\nu)\lesssim  a^{1+1/p}\norm {U_\mu-U_\nu}_{L^p(B_{K+1/a}(0))}+\sup_{R\geq 0,z_0\in\IC}\nu\big(\overline{B_{R+1/a}(z_0)\setminus B_R(z_0)}\big)
 \end{align*} 
\end{theorem}

In the same manner it is possible to show an analogue for the classical Kolmogorov distance between 2-dimensional distribution functions, see Corollary \ref{cor:SmoothingKolm}. For measures $\mu, \nu$ on $\IR$, where $\nu$ has a bounded density, Dinh and Vu \cite{D17} showed another direct relation of similar type
\begin{align*}
 \abs{\mu (I)-\nu (I) }\lesssim \norm{U_\mu -U_\nu}_{L^\infty(\supp \nu)}^{1/2}
\end{align*}
for all intervals $I\subseteq \IR$ and it was used to show a rate of convergence in Wigner's Semicircular Law and the Marchenko-Pastur Law. Theorem \ref{prop:Smoothing} may be of independent interest, since it can be considered as a complex counterpart of other smoothing inequalities of distributions $\mu,\nu$ on the real line. For instance in the case of Fourier transforms $\phi_\mu(t)=\int e^{itx}d\mu(x)$, the well known \emph{Berry-Essen inequality}
\begin{align}\label{eq:BerryEsseen}
  \sup_{x\in\IR}\abs{(\mu-\nu)((-\infty,x]) }\lesssim \int_{-a}^a \abs{\frac{\phi_\mu(t)-\phi_\nu(t)}{t}}dt+\sup_{x\in\IR}\nu((x,x+c/a])
 \end{align}
leads to a rate of convergence of order $1/\sqrt n$ in the Central Limit Theorem, when choosing $\nu=\mathcal N (0,1)$ and $\mu=\IP_{S_n}$ for the normalized sum $S_n=n^{-1/2}\sum_{k=1}^n X_k$ of i.i.d. random variables $X_k$ with $\IE X_1 =0, \IE X_1^2 =1$ and finite third moment $\IE X_1^3<\infty$. In Random Matrix Theory, \emph{Bai's inequality} is a handy tool to profit from control of Stieltjes' transforms $m_\mu(z)=\int\frac 1 {x-z} d\mu(x)$ that can be simplified to
 \begin{align}\label{eq:BaisIneq}
  \sup_{x\in\IR}\abs{(\mu-\nu)((-\infty,x]) }\lesssim \int \abs{m_\mu-m_\nu}(t+i/a)dt+\sup_{x\in\IR}\nu((x,x+c/a]). 
 \end{align}
Roughly speaking, \cite{BS10spectral} uses $a\simeq \sqrt n$ to show a rate of convergence of order $1/\sqrt n$ for the Kolmogorov distance in Wigner's semicircle law under finite sixth moment condition. Using an improved, but more involved smoothing inequality, it is shown in \cite{GT16Optimal} that the optimal rate of convergence to the semicircle distribution is given by $\mathcal O (1/n)$, even in the interior of the bulk. To our knowledge, the best rate of convergence of the non-averaged ESD to the Semicircle Law is given by $\mathcal O (\log^2(n)/n)$ obtained in \cite[Equation (1.10)]{GNTT}.

\subsection{Rate of Convergence}\label{subsec:roc} All smoothing inequalities \eqref{eq:BerryEsseen}, \eqref{eq:BaisIneq} and Theorem \ref{prop:Smoothing} are used to derive convergence rates under moment conditions and they share the essential structure of bounding the Kolmogorov distance by the distance of certain integral-transforms and an additional maximal annulus probability of width $\mathcal O (1/a)$ with respect to the ``limit distribution''. Regarding Theorem \ref{prop:Smoothing}, we consider the distributions $\mu=\mu_n, \nu=\mu_\infty$ from the introduction and choose $a=\sqrt n, K=1$. In this case we see that the remainder term is of order $n^{-1/2}$ and a rate of convergence for $D_n$ follows.

It is important to carefully distinguish between events holding with overwhelming probability uniformly in $z$ and uniform events that hold w.o.p.. The former leads to Local Circular Laws like \cite[Theorem 20]{TV15uni} (see also Corollary \ref{cor:localcirc} below) and hence do not imply the latter, which is an estimate on $D_n$. Contrary to Local Circular Laws, a bound on $D_n$ w.o.p. allows to choose the (worst) ball $B_R(z_0)$ depending on the random sample of the eigenvalues $(\lambda_j (X(\omega)/\sqrt n))_j$, cf. Figure \ref{fig:Samples}.
Similarly, the statement of Proposition \ref{prop:ConcLogPotAEK} should not be confused with an assertion about the uniform term $\sup_{z\in B_K(0)}\abs{U_n(z)-U_\infty (z) }$, since it equals $\infty$ whenever an eigenvalue lies in $B_K(0)$. Due to this fact one cannot simply take $p=\infty$ in Proposition \ref{prop:Smoothing} in order to obtain the following result.

\begin{theorem}\label{thm:rateofconv}
If condition (A) holds, then for every (small) $\epsilon>0$ and (large) $Q>0$ 
\begin{align}\label{eq:rateofconv}
\IP(D_n\leq n^{-1/2+\epsilon})\geq 1-n^{-Q}
\end{align}
holds for sufficiently large $n$, where $\displaystyle D_n=\sup_{z_0\in\IC,R>0}\abs{(\mu_n-\mu_\infty) (B_R(z_0)) }$.
\end{theorem}

By virtue of Corollary \ref{cor:SmoothingKolm}, the following Kolmogorov distance analogue holds.
\begin{theorem}\label{thm:rateofconvKolm}
If condition (A) holds, then for every $\epsilon,Q>0$ 
\begin{align}\label{eq:rateofconvKolm}
\IP( d_n\leq n^{-1/2+\epsilon})\geq 1-n^{-Q}
\end{align}
holds for sufficiently large $n$, where $\displaystyle d_n=\sup_{s,t\in\IR}\abs{(\mu_n-\mu_\infty)((-\infty,s]+i(-\infty,t]) }$.
\end{theorem}

Invoking Proposition \ref{thm:ConcLogPot}, we prove a rate of convergence result weakening the conditions of the last statements at the cost of excluding sets close to the edge.

\begin{theorem}\label{thm:rateofconvBulk}
If condition (B) holds, then for every $\epsilon,\tau,Q>0$ 
\begin{align*}
\IP(D_n ^{\circ}\leq n^{-1/2+\epsilon})\geq 1-n^{-Q}
\end{align*}
holds for sufficiently large $n$, where $\displaystyle D_n ^{\circ}=\sup_{B_R(z_0)\subseteq B_{1-\tau}(0)}\abs{(\mu_n-\mu_\infty)(B_R(z_0))}$.
\end{theorem}

A similar result also holds true for products of independent matrices. In this note however we stick to the limiting Circular Law, since new methods are needed to derive the optimal rate given by products of Ginibre Matrices that shall be discussed in the subsequent work \cite{Jalowy}.

\subsection{Discussion} 

Tao and Vu \cite{TV08Circular} showed that with probability 1 the Kolmogorov distance $ d_n$ of the 2-dimensional distribution functions is of order $n^{-\eta}$ for some small, unknown $\eta>0$, which holds for finite $2+\epsilon$-moments of the entries. Comparing this to Theorem \ref{thm:rateofconvKolm}, we see that a nearly optimal rate of convergence is obtained in \eqref{eq:rateofconvKolm} which holds with overwhelming probability. On the other hand a stronger moment assumption for the entries is needed. In particular, this explicit rate of convergence gives a partial answer to an open problem mentioned in \cite{TV09Littlewood}. As already discussed above, non-uniform rates can be read off from Local Circular Laws \cite{BYY14Local,BYY14LocalEdge,TV15uni,GNT17Local} and fluctuations of linear spectral statistics \cite{RS06LinStat,KOV18}. Note that these results deal with certain classes of smooth functions, whereas the metric $D$ uniformly covers classes of non-smooth indicator functions. One may ask for other function classes, i.e. Lipschitz functions corresponding to rates of convergence in terms of Wasserstein distances or bounded Lipschitz metric. However Lipschitz functions may have uncontrollable Laplacians, which is essential for the logarithmic potential approach due to \eqref{eq:logPot}.

In the special case of Gaussian entries, i.e. for the Ginibre ensemble, pointwise convergence of the density of $\bar\mu_n$ has been also discussed in \cite{AC18, TV15uni}, similar to the integrated version \eqref{eq:GinibreRateBulk}. Furthermore $\IP$-a.s. convergence rates of order $\sqrt{\log n }/n^{1/4}$ in $p$-Wasserstein distance for $1\leq p\leq 2$ have been proven in \cite{MM14rate}. Recently this has been extended by O'Rourke and Williams \cite{OW19} to matrices satisfying a moment matching condition, where a non-optimal rate of $\mathcal O (n^{-1/4+o(1)})$ in $1$-Wasserstein distance has been shown. Though the Wasserstein distance is not directly comparable to $D$, both optimal rates are expected to be $n^{-1/2}$ up to logarithmic factors. Moreover, for Ginibre matrices and centered balls, the fluctuation around the (deterministic) rate has been studied in \cite{FL20} and for a restricted class functions with controlled Laplacian, the rate can be improved to $1/n$, see \cite{Lambert}.

More generally, Chafa\"\i, Hardy and Ma\"\i da studied invariant $\beta$-ensembles with external potential $V$ instead of independent-entry matrices in \cite{CHM16}. Their result implies a rate of convergence to the limiting measure with density $c\Delta V$ of order $\sqrt{ \log n / n}$ with respect to the bounded Lipschitz metric and the $1$-Wasserstein distance.
The paper \cite{CHM16} is also based on an inequality between distances of measures to their energy, i.e. integrated logarithmic potential, similar to Theorem \ref{prop:Smoothing}. 
However it relies critically on the existence of a confining potential, hence a joint probability density function for the eigenvalues. Note that their result is given for a Coulomb gas point process in arbitrary dimension $d>1$, yielding a bound of order $n^{-1/d}$ up to logarithmic factors. This coincides with the rate of order $1/n$ for the semicircle law for $d=1$ as well as the optimal order $1/\sqrt{n}$ in the Circular Law and can also be interpreted as mentioned in the introduction. Similar questions in this context of $\log$-gases, but for the non-uniform variant of $D$ (the discrepancy) have been addressed in \cite{Serfaty}.

\subsection{Application to Random Polynomials}\label{sec:Polynomials} 
The Smoothing Inequality can also be applied to the empirical distribution of roots of random polynomials in order to obtain the same rate of convergence to the Circular Law as before. In the previous sections we considered the roots of the characteristic polynomial of a random matrix, where the coefficients of the polynomial exhibit specific dependencies. We begin by replacing the independence condition on the matrix entries by independent coefficients in the polynomial.
\begin{definition}\label{defi:RandomPolynomials}
Given $n\in\IN$ many complex numbers $c_0,\dots,c_n$ and i.i.d. centered complex random variables $\xi_0,\dots,\xi_n$ with $\IE \abs{\xi_{k}}^2=1$ and $\IE \abs{\xi_{0}}^{2+\delta}<\infty$ for some $\delta>0$, we define the random polynomial $f_n:\IC\to\IC$ by
\begin{align*}
f_n(z)=\sum_{k=0}^{n} c_k \xi_k z^k.
\end{align*}
In particular we will work with so called \emph{Weyl (or Flat) polynomials} $f_n^W$ corresponding to $c_k=\sqrt{ n^k/{k!}}$. By analogy to the Introduction, we associate to a random polynomial $f_n$ its multiset of zeros $\Lambda:=\{\lambda\in\IC:f_n(\lambda)=0\}$ taking their multiplicities into account and its empirical measure given by 
\begin{align*}
\mu_{f_n}=\frac 1 n\sum_{\lambda\in \Lambda}\delta_{\lambda}.
\end{align*}
\end{definition}

It should be remarked that $\mu_{f_n}$ is not necessarily normalized, since a random polynomial may have degree deg$(f_n)<n$. Unsurprisingly this does not affect the large $n$ limit, since $n-$deg$(f_n)\in\mathcal O (1)$ $\IP$-a.s.. As in \cite{IZ13uni}, we may always assume $\IP(\xi_0=0)=0$, since otherwise we may restrict ourselves to the events $\{$deg$(f_n)=k, \min\{j\leq n :\xi_j\neq 0\}=l \}$ with fixed degree and fixed amount of zero roots.

The Circular Law for the empirical measure of the roots of Weyl polynomials has been established in \cite{KZ14Polyn} by Kabluchko and Zaporozhets, see also \cite{FH99} for the Gaussian case, stating
\begin{align*}
\mu_{f_n^W}\Rightarrow\mu_\infty\quad \IP\text{-a.s.}.
\end{align*}
Note that their result holds for much more general random analytic functions and under the much weaker condition of the coefficients having finite logarithmic moments $\IE\log(1+\abs{\xi_0})<\infty$. 

We aim to quantify this result by showing a rate of these convergences of order $n^{-1/2+\epsilon}$ by using results about logarithmic potentials. Since local universality for certain random polynomials has been proven in by Tao and Vu using concentration of logarithmic magnitudes $\log \abs{f_n}$, we can apply the same methods as before. We denote $U_n=-\frac 1 n \log\abs{f_n}$ and rephrase \cite[Lemma 12.1]{TV15Polynomials}: For every $\epsilon,\delta,\tau,Q>0$ there exist a constant $c>0$ such that 
\begin{align}\label{eq:ConcLogPotPoly}
\IP\left(\abs{U^W_n(z)-U_\infty(z)+1/2}\leq c n^{-(1-\epsilon)}\right)\geq 1-n^{-Q}
\end{align}
holds uniformly for $n^{-1/2+\delta}\le\abs z\le 1+\tau$. The origin has to be avoided, since the distribution of $U_n^W(0)=-\frac 1 n\log\abs{\xi_0}$ around $0$ is not necessarily converging. In particular, the bound \eqref{eq:ConcLogPotPoly} will not hold in $z=0$ if $\IP(\xi_0=0)>0$. Due to the application of the Monte Carlo method we still need a technical assumption on the concentration of $\xi_0$ near $z=0$ in the following rate of convergence result which we deduce from a variant of the smoothing inequality, Theorem \ref{prop:Smoothing}.

 \begin{theorem}\label{thm:rateofconvPolynom}
If $\IE\abs{1/\xi_0}^\delta<\infty$ for some $\delta>0$, then for every $\epsilon,Q>0$ and sufficiently large $n$ we have 
\begin{align*}
\IP(D (\mu_n^W,\mu_\infty)\leq n^{-1/2+\epsilon})\geq 1-n^{-Q}.
\end{align*}
\end{theorem}

It seems likely that other polynomials, like elliptic polynomials, converge at the same rate to their corresponding limit root distributions, but we focus on Circular Laws in this work.

\section{Proofs of the Smoothing Inequalities}\label{sec:Smoothing}

We will proof the following slightly more general statement that covers all variants of the smoothing inequality which we will need to establish Theorems \ref{thm:rateofconv}-\ref{thm:rateofconvPolynom}.
 \begin{theorem}\label{thm:SmoothingThm}
Let $\mu,\nu$ be probability measures on $\IC$ with logarithmic potentials $U_\mu, U_\nu$ respectively (i.e. the distributional Poisson equation of \eqref{eq:logPot} holds). Fix $1\le p\leq \infty$, $z^*\in\IC$, $K>0$ and $\eta\ge 0$. For any $a>1$, define the annuli $V=B_{K}(z^*)\setminus B_{2\eta/a}(z^*)$ and $V '=B_{K+2/a}(z^*)\setminus B_{\eta/a}(z^*)$, s.t. it holds
\begin{align*}
 D(\mu,\nu) \lesssim&  a^{1+1/p}\norm {U_\mu-U_\nu}_{L^p(V')} +\mu(V^c)+\nu(V^c)\\
 &+\sup_{R\geq 0,z_0\in\IC}\nu\left(z\in V': R\leq \abs{z-z_0 }\leq R+\max(2,\eta)/a\right).
 \end{align*}  
 \end{theorem}
Here, $\eta\neq 0$ is only needed for the applications to random polynomials, where the logarithmic potential near the origin cannot be controlled. We retrieve Theorem \ref{prop:Smoothing} by taking $\eta=0$, $z^*=0$, $\nu(V^c)=0$, replacing $a$ by $2a$ for simplicity and noting that for probability distributions we estimate
\begin{align*}
\mu(V^c)=(\nu-\mu)(V)\le\sup_{R\ge 0, z_0\in\IC}\abs{ (\mu-\nu)(B_R(z_0)\cap V)}
\end{align*}
in the following proof, see \eqref{eq:start} and \eqref{eq:Zwischenabschaetzung} below.
 
\begin{proof}
First, note that 
\begin{align}\label{eq:start}
 \sup_{R\ge 0, z_0\in\IC}\abs{ (\mu-\nu)(B_R(z_0))} \le \sup_{R\ge 0, z_0\in\IC}\abs{ (\mu-\nu)(B_R(z_0)\cap V)} +\mu(V^c)+\nu(V^c),
 \end{align} 
hence we have to estimate the first term. Fix some $a>1$, let $\phi\in\mathcal{C}^\infty(\IR)$ be nonnegative with $\supp \phi \subseteq [-1,1]$ and $\int\phi=1$, and define $\phi_a(\rho)=a\phi(a\rho)$. For arbitrary $R>0$ and $z_0\in\IC$ we mollify the indicator function appearing in $D(\mu,\nu)$ via the rotationally invariant approximation
\begin{align*}
 f_1(z):&=\left(\eins_{(-\infty, R-1/a]}\ast\phi_{a}\right)(\abs{z-z_0} )\\
 &\leq\eins_{B_R(z_0)}(z)\\
 &\leq\left(\eins_{(-\infty,R+1/a]}\ast\phi_a\right)(\abs{z-z_0} )=:f_2( z),
\end{align*}
where we choose $f_1\equiv 0$ if $R\le 2/a$ for smoothness reasons. Furthermore we will approximate $\eins_V$ by smooth functions $h_1$ from inside and by $h_2$ from outside, more precisely define
\begin{align*} 
h_1(z):=\begin{cases}
\left( (\eins_{[5\eta/2a,\infty)}\ast\phi_{2a/\eta}) \cdot( \eins_{(-\infty,K-1/a]}\ast\phi_{a})\right)(\abs{z-z^*} )&\text{, if } \eta>0,\\
\eins_{(-\infty,K-1/a]}\ast\phi_{a}(\abs{z-z^*} )&\text{, if } \eta=0 ,         
         \end{cases}\\
 h_2(z):=\begin{cases}
\left( (\eins_{[3\eta/2a,\infty)}\ast\phi_{2a/\eta}) \cdot( \eins_{(-\infty,K+1/a]}\ast\phi_{a})\right)(\abs{z-z^*} )&\text{, if } \eta>0,\\
\eins_{(-\infty,K+1/a]}\ast\phi_{a}(\abs{z-z^*} )&\text{, if } \eta=0 .         
         \end{cases}
\end{align*}
We apply $h_1f_1\leq \eins_{B_R(z_0)\cap V}$ and integration by parts (in other words we use the definition of the distributional Poisson equation \eqref{eq:logPot}) back and forth to obtain
\begin{align*}
&\mu(B_R(z_0)\cap V)\geq\int h_1f_1d\mu=-\frac{1}{2\pi}\int\Delta (h_1f_1) U_\mu d\lebesgue\\
&=-\frac{1}{2\pi}\int\Delta (h_1f_1) (U_\mu-U_\nu) d\lebesgue-\int(\eins_{B_R(z_0)\cap V}-h_1f_1)d\nu+\int \eins_{B_R(z_0)\cap V}d\nu,
\end{align*}
where $\lebesgue$ denotes the Lebesgue measure on $\IC$. A rough estimate of the error of approximation yields for the second term 
\begin{align}
 \int(\eins_{B_R(z_0)\cap V}-h_1f_1)d\nu&\leq  \nu\left(z\in V': R-2/a\leq \abs{z-z_0 }\leq R\right) +\nu( V'\setminus V) \label{eq:remainderV}\\
 &\le3\sup_{R\geq 0,z_0\in\IC}\nu\left(z\in V': R\leq \abs{z-z_0 }\leq R+\max(2,\eta)/a\right)\nonumber\\
 &=:3M_\nu(a).\nonumber
\end{align}
We use H\"older's inequality to estimate the first term, implying
\begin{align}\label{eq:Losing}
(\mu-\nu)(B_R(z_0)\cap V)\geq -\frac{1}{2\pi}\norm{\Delta (h_1f_1) }_{L^q}\norm{U_\mu-U_\nu}_{L^p}-3M_\nu(a),
\end{align}
where $L^p=L^p(V')$, $L^q=L^q(V')$ (we omit $V'$ in the sequel), $1/p+1/q=1$ and $R>0,z_0\in\IC$ are still arbitrary. Noting $\mu(B_R(z_0)\cap V)\leq\int h_2f_2d\mu$ and taking the same route for $h_2f_2$ as for $h_1f_1$, we obtain the same upper bound, i.e.
\begin{align}
-&\frac{1}{2\pi}\norm{\Delta (h_1 f_1) }_{L^q}\norm{U_\mu-U_\nu} _{L^p}-3M_\nu(a) \nonumber\\
\leq& (\mu-\nu)(B_R(z_0)\cap V)\label{eq:Zwischenabschaetzung} \\
\leq& \frac{1}{2\pi}\norm{\Delta (h_2 f_2) }_{L^q}\norm{U_\mu-U_\nu} _{L^p}+3M_\nu(a).\nonumber
\end{align}
Therefore it remains to control 
\begin{align*}
\norm{\Delta (h_jf_j) }_{L^q}\leq \norm{h_j\Delta f_j }_{L^q}+2\norm{\nabla h_j \cdot\nabla f_j }_{L^q}+\norm{f_j\Delta h_j }_{L^q}.
\end{align*}
We see that the supports of all three functions are (unions of) annulus-segments, e.g. $V'\cap (B_{R+2/a}(z_0)\setminus B_R(z_0) )$ for $h_2\Delta f_2$, with length at most $2\pi(K+2/a)$ and the width equal to $\max(2,\eta)/a$. Hence uniformly in $R>0$ and $z_0\in \IC$, the size of the area of integration is bounded by $cK\max(2,\eta)/a$ and we arrive at
\begin{align*}
 \norm{\Delta (h_jf_j) }_{L^q}\leq (cK\max(2,\eta)/a)^{1/q}\left(\norm{\Delta f_j }_{L^\infty}+2\norm{\nabla h_j \cdot\nabla f_j }_{L^\infty}+\norm{\Delta h_j }_{L^\infty}\right).
\end{align*}
With our choice of $f_j$ and $h_j$, the radial derivatives become fairly simple, e.g.
\begin{align*}
\partial_r (f_2 (z+z_0))=\partial_r\int_{\abs z -R-1/a}^\infty\phi_a(\rho)d\rho=-a\phi(a\abs z - aR-1).
\end{align*}
Due to the rotational symmetry of $f_2$, we have $\norm{\nabla f_2}_{L^\infty}\leq \norm{\phi_a}_{L^\infty}\lesssim a$ and again exploiting rotational symmetry it follows that the maximal curvature is attained in radial direction, i.e. 
\begin{align}\label{eq:Laplacefunendlich}
\norm{\Delta f_2 }_{L^\infty}=\sup_{r>0}\abs{\partial_r^2 f_2 (z_0+r) }=a^2\norm{\phi' }_{L^\infty} . 
\end{align} 
The same bounds hold for $j=1$ and $h_j$, where for $h_j$ we replace $a$ by $a\max(1,2/\eta)$ if $\eta>0$.
Finally we conclude
\begin{align}
 \norm{\Delta h_jf_j }_{L^q}\lesssim  a^2(K/a)^{1/q}\lesssim K^{1-1/p}a^{1+1/p},\label{eq:Laplacef2}
\end{align}
where the implicit constant in the last $\lesssim$ depends on $p,\eta$ and $\phi$ only. More precisely, one may combine all previous estimates and choose $\phi$ with sufficiently small derivative such that the constant can be chosen to be $c=200$ for $\eta=0$ or $c=800\eta^{-2}$ for $\eta<2$. The claim now follows from taking the supremum over $R>0$ and $z_0\in\IC$ in \eqref{eq:Zwischenabschaetzung}.
\end{proof}

In fact if we restrict ourselves to a certain region, we obtain a local smoothing inequality that makes it possible to invoke Proposition \ref{thm:ConcLogPot}.

\begin{corollary}\label{cor:SmoothingLocal}
Let $\mu,\nu$ be probability measures on $\IC$ with logarithmic potentials $U_\mu, U_\nu$ respectively, and fix some $z^*\in\IC$, $K,\tau>0$ and $1\le p\leq \infty$. There exists a constant $c>0$ such that for any $a>1\wedge \tau^{-1}$
\begin{align*}
 \sup_{B_R(z_0)\subseteq B_{K-\tau}(z^*)}\abs{ (\mu-\nu)(B_R(z_0))} \le c  a^{1+1/p}\norm {U_\mu-U_\nu}_{L^p(B_{K}(z^*))} \\
 +\sup_{R\geq 0,z_0\in\IC}\nu\left(z\in B_{K}(z^*): R\leq \abs{z-z_0 }\leq R+2/a\right).
 \end{align*}  
\end{corollary}
\begin{proof}
 Replace $K$ by $K-\tau$, set $\eta=0$ and note that the cutoff $h$ in the previous proof is not necessary anymore, since $B_R(z_0)\subseteq B_{K-\tau}(z^*)$. In particular the last term of \eqref{eq:remainderV} does not exist and we only need to bound \eqref{eq:Laplacefunendlich}.
\end{proof}

Albeit we will only use this inequality for $K=1, z^\ast =0$, this inequality shows that the local distance of the measures only depends on the local distance of the logarithmic potentials. Girko's Hermitization Trick however transforms it to a highly nonlocal problem, taking the whole support (or spectrum for $\mu=\mu_n$) into account.

Moreover, the method of proof extends to the case of the classical Kolmogorov distance between 2-dimensional distribution functions.
 
\begin{corollary}\label{cor:SmoothingKolm}
Let $\mu,\nu$ be probability measures on $\IC$ with $\supp\nu\subseteq [-K,K]^2$ for some $K>0$, let $U_\mu, U_\nu$ be their logarithmic potentials and fix some $\tau>0$ and $1\le p\leq \infty$. There exists a constant $c>0$ such that for any $a>1$
\begin{align*}
 \sup_{s,t\in\IR}\abs{(\mu-\nu)((-\infty,s]+i(-\infty,t]) }&\leq c a^{1+1/p}\norm {U_\mu-U_\nu}_{L^p([-K-\tau,K+\tau]^2)} \\
 &+3\sup_{s,t\in\IR}\nu(([s,s+2/a]+i \IR)\cup (\IR+i[t,t+2/a])).
 \end{align*} 
\end{corollary}

\begin{proof}
We exploit the same ideas from the proof of Theorem \ref{thm:SmoothingThm} by finding a substitute of \eqref{eq:Zwischenabschaetzung}. We continue with the same notation, where $[-K,K]^2$ takes the role of $V$ and $\tau$ corresponds to $2/a$. Define now
\begin{align*}
 f_1(z):&=\eins_{(-\infty, s-1/a]}\ast\phi_{a}(\Re z) \cdot\eins_{(-\infty, t-1/a]}\ast\phi_{a}(\Im z)\\
 &\leq\eins_{(-\infty,s]+i (-\infty,t]}(z)\\
 &\le\eins_{(-\infty, s+1/a]}\ast\phi_{a}(\Re z)\cdot \eins_{(-\infty, t+1/a]}\ast\phi_{a}(\Im z)=:f_2(z),
\end{align*}
and $h(z)=\eins_{[-K-\tau/2,K+\tau/2]}\ast\phi_{\tau/2}(\Re z) \cdot\eins_{[-K-\tau/2,K+\tau/2]}\ast\phi_{\tau/2}(\Im z)$. Here, since $\nu$ has compact support, we do not need $h_1$ in order to restrict ourselves to $V$. By similar arguments as above, e.g. $hf_1\le\eins_{(-\infty,s]+i (-\infty,t]}$, we obtain
\begin{align*}
(\mu-\nu)((-\infty,s]+i(-\infty,t])\ge -\frac{1}{2\pi}\norm{\Delta (hf_1) }_{L^q}\norm{U_\mu-U_\nu}_{L^p}-M_\nu(a),
\end{align*}
where now $M_\nu(a)=\sup_{s,t\in\IR}\nu(([s,s+2/a]+i \IR)\cup (\IR+i[t,t+2/a]))$ and we abbreviated $L^p=L^p([-K-\tau,K+\tau]^2)$, $L^q=L^q([-K-\tau,K+\tau]^2)$. For a short moment, consider 
\begin{align*}
f_1^0(z)=(\eins_{[-K+1/a,K-1/a]}\ast\phi_{a}(\Re z))\cdot (\eins_{[-K+1/a,K-1/a]}\ast\phi_{a}(\Im z))
\end{align*}
which analogously to the idea mentioned before Corollary \ref{cor:SmoothingLocal} yields
\begin{align*}
 1-\mu([-K,K]^2)=(\nu-\mu)([-K,K]^2)\le\frac{1}{2\pi}\norm{\Delta (f_1^0) }_{L^q}\norm{U_\mu-U_\nu}_{L^p}+2M_\nu(a).
\end{align*}
We conclude
\begin{align*}
-&\frac{1}{2\pi}\norm{\Delta (h f_1) }_{L^q}\norm{U_\mu-U_\nu} _{L^p}-M_\nu(a) \\
\leq& (\mu-\nu)((-\infty,s]+i (-\infty,t]) \\
\leq& \frac{1}{2\pi}\left(\norm{\Delta (h f_2) }_{L^q}+\norm{\Delta (f^0_1) }_{L^q}\right)\norm{U_\mu-U_\nu} _{L^p}+3M_\nu(a).
\end{align*}
Consequently it remains to derive similar estimates $\norm{\Delta (hf_j)}_{L^q}\lesssim a^{1+1/p}$ using the same arguments as before. We omit the details here.
 \end{proof}

\section{Proofs of the Rates of Convergence}\label{sec:Rates}

\begin{proof}[Proof of Theorem \ref{thm:rateofconv}]
Without loss of generality $\epsilon<4$, we choose $p>4/\epsilon$ and apply Theorem \ref{prop:Smoothing} to $\mu=\mu_n,\nu=\mu_\infty, K=1$ and $a=\sqrt n$, 
 \begin{align*}
 D_n\lesssim n^{1/2+\epsilon/2}\norm {U_n-U_\infty}_{L^p(B_{1+\tau}(0))}+\sup_{R\geq 0,z_0\in\IC}\mu_\infty\left( R\leq \abs{\cdot-z_0 }\leq R+2n^{-1/2}\right).
 \end{align*}
Since $\mu_\infty$ has bounded support and bounded density it is clear that the second term is of order $\mathcal O (n^{-1/2})$. In order to obtain a bound of the $L^p(B_{1+\tau}(0))$-norm of the log potentials from the pointwise estimate in Proposition \ref{prop:ConcLogPotAEK}, we adapt the Monte Carlo sampling method which was used in \cite{TV15uni} (in a different form); we set $I(z)= \abs{U_n(z)-U_\infty(z)}$ and approximate 
\begin{align*}
\fint I(z)^pdz:=\frac 1{\pi(1+\tau)^2}\int_{B_{1+\tau}(0)}I(z)^pdz \approx\frac{1}{m}\sum_{j=1}^m I(z_j)^p=:S_m,
\end{align*}
where $(z_j)_{j=1,\dots,m}$ are independent random variables (also independent of $X_{ij}$) uniformly distributed on $B_{1+\tau}(0)$. More precisely we will show that for every $Q>0$
\begin{align}\label{eq:bzz1}
\abs{\fint I(z)^pdz-S_m}^{1/p}\lesssim n^{-1}
\end{align}
as well as
\begin{align}\label{eq:bzz2}
\abs{S_m}^{1/p}\lesssim n^{-1+\epsilon/2}
\end{align}
holds with probability at least $1-n^{-Q}$ for some large $n$-dependent $m$. Assuming \eqref{eq:bzz1} and \eqref{eq:bzz2} are true, there exist constants $c_1,c_2,c_3,c_4>0$ such that
\begin{align*}
 &\IP(D_n\geq c_1n^{-1/2+\epsilon})\\
 &\leq \IP\left( c_2n^{1/2+\epsilon/2}\Big(\Big|\fint I(z)^pdz-S_m\Big|^{1/p}+\abs{S_m}^{1/p}\Big)+ c_2n^{-1/2}\geq  c_1n^{-1/2+\epsilon}\right)\\
 &\leq \IP\Big( \Big|\fint I(z)^pdz-S_m\Big|^{1/p}\geq c_3n^{-1+\epsilon/2}\Big) + \IP\Big(\abs{S_m}^{1/p}\geq c_4n^{-1+\epsilon/2}\Big) \\
 &\leq n^{-Q}
\end{align*}
proving the claim.\\
Let us turn to the proof of \eqref{eq:bzz1}. First, we restrict ourselves to the set of polynomially bounded eigenvalues. On the one hand the largest absolute value of eigenvalues $\abs\lambda_{\max } $ is bounded by the largest singular value $s_{\max } $ and on the other hand for every $Q>0$ we have
\begin{align}\label{eq:s_max}
 \IP(s_{\max }\geq n^{(Q+1)/2})\leq \frac 1 {n^{Q+1}}\IE\norm{X/\sqrt n}^2\leq \frac 1 {n^{Q+2}} \sum_{ij}^n\IE\abs{X_{ij}}^2\leq n^{-Q},
\end{align}
where the operator norm $\norm\cdot$ has been estimated by the Hilbert Schmidt norm. We freeze the coefficients $X_{ij}$ and use Chebyshev's inequality for the probability measure conditioned on $X$
\begin{align*}
 \IP\Big( \Big|S_m-\fint I(z)^pdz\Big|^{1/p}\geq \frac c n\Big|X\Big)\leq \frac{n^{2p}}{c^{2p}} \Var (S_m|X)\leq\frac{n^{2p}}{mc^{2p}} \Var (I(z_1)^p|X).
\end{align*}
The variance of $I(z_1)^p$ given $X$ is given by
\begin{align*}
 \Var (I(z_1)^p|X)\leq \IE(I(z_1)^{2p}|X)\leq \fint \abs{U_n(z)}^{2p}+\abs{U_\infty(z)}^{2p}dz.
\end{align*}
If we assume the eigenvalues $\lambda_1,\dots,\lambda_n$ to be fixed and use Jensen's inequality, we may estimate
\begin{align*}
 \fint\abs{U_n(z)}^{2p}dz&\leq \frac 1 n \sum_{j=1}^n\fint\abs{\log\abs{\lambda_j-z}}^{2p}dz\\
 &\leq  \frac c n \sum_{j=1}^n \iint_{B_{1+\tau}(-\lambda_j)}r\abs{ \log r}^{2p}drd\phi\\
 &\leq c_p (1+\tau+\abs\lambda_{\max})\log^{2p}(1+\tau+\abs\lambda_{\max})\\
 &\leq c_pn^{(Q+1)/2}\log^{2p}n
\end{align*}
for some $p$-dependent constant $c_p$. Similarly we get $\fint\abs{U_\infty(z)}^{2p}dz=c_p$. Now choose $m:=n^{2p+3Q/2+1}$ and putting the estimates together we have shown
\begin{align*}
 \IP&\Big(\Big|\fint I(z)^pdz-S_m\Big|^{1/p}\geq cn^{-1}\Big)\\
 &\leq  \IE\Big(\IP\Big(\Big\{\Big|\fint I(z)^pdz-S_m\Big|^{1/p}\geq \frac cn\Big\}\cap\Big\{\abs\lambda_{\max }\leq n^{\frac{Q+1}{2} }\Big\}\Big|X\Big)\Big)+n^{-Q}\\
 &\leq c n^{-Q}.
\end{align*}
It remains to show \eqref{eq:bzz2}. To this end we use Proposition \ref{prop:ConcLogPotAEK} with an adjusted error probability stating 
\begin{align}\label{eq:poflogorder}
 \IP(I(z)\geq cn^{-1+\epsilon/2})\leq n^{-2p-5Q/2-1}
\end{align}
uniformly in $B_{1+\tau}(0)$. If $I(z_j)\leq n^{-1+\epsilon/2}$ for all $j=1,\dots,n$ then $\abs{S_m}^{1/p}\leq n^{-1+\epsilon/2}$ which implies 
\begin{align*}
\IP( \abs{S_m}^{1/p}\geq n^{-1+\epsilon/2}) &\le\sum_{j=1}^m\IP(I(z_j)\geq n^{-1+\epsilon/2})\\
&\leq cmn^{-2p-5Q/2-1}=cn^{-Q}.
\end{align*}
The proof is now complete, since these constants may be absorbed by the $n^{-Q}$ (respectively $n^\epsilon$-)term for some slightly larger $Q$ (respectively smaller $\epsilon$).
\end{proof}

Analogously, Theorem \ref{thm:rateofconvKolm} follows from Corollary \ref{cor:SmoothingKolm} and Theorem \ref{thm:rateofconvBulk} follows from Corollary \ref{cor:SmoothingLocal}. The details are exactly the same as above and we skip them. Moreover using the same techniques it is possible to show the following version of a Local Circular Law. Compared to \cite{GNT17Local} it improves the statement to hold with overwhelming probability but replaces the constant $\norm{\Delta f}_{L^1}$ by $\norm{\Delta f}_{L^q}$ and is stated for a single matrix, instead for a product of $m$ many.

\begin{corollary}[Local Circular Law]\label{cor:localcirc}
Let $q>1$, $z_0\in B_{1+\tau^{-1}}(0)$ with $\abs{1-\abs {z_0} }\ge \tau$, $f:\IC\to\IR_+$ be a bounded smooth function, which is compactly supported with $\norm{f'}_{L^\infty}\leq n^c$ for some constant $c>0$. Define the function $f_{z_0}(z):=n^{2s}f((z-z_0)n^s)$ which zooms into $z_0$ at speed $s\in(0,1/2)$. For any $Q>0$ there exist a constant $c>0$ such that 
\begin{align*}
 \IP\left(\abs{\frac 1 n \sum_{j=1}^nf_{z_0}(\lambda_j)-\int_{\IC}f_{z_0}(z)d\mu_\infty(z)}	\leq \frac{c\log ^4 n}{n^{1-2s}}\norm{\Delta f}_{L^q} \right)\geq 1-n^{-Q}.
\end{align*}
\end{corollary}

Recalling the discussion in section \ref{sec:Main}, $z_0$ and  $f$ are not allowed to depend on $\omega$ here.

\begin{proof}
As in the proof of Theorem \ref{prop:Smoothing}, integration by parts yields
\begin{align*}
 \frac 1 n \sum_{j=1}^nf_{z_0}(\lambda_j)-\int_{\IC}f_{z_0}(z)d\mu_\infty(z)=-\frac{n^{2s}}{2\pi}\int_{\IC}\Delta f(z) \left(U_n(z)-U_\infty(z)\right)dz.
\end{align*}
After applying H\"older's inequality as was done in \eqref{eq:Losing}, it remains to show the estimate $\norm{U_n-U_\infty }_{L^p}\lesssim \log^4 n/n$ similar to the proof of Theorem \ref{thm:rateofconv}. More precisely, the error \eqref{eq:bzz1} of the Monte Carlo sampling is already sufficiently small and the estimate of \eqref{eq:bzz2} can be improved by using Proposition \ref{thm:ConcLogPot} in the analogue of \eqref{eq:poflogorder}.
\end{proof}

We now turn to an application for random polynomials. The proof does not differ much from those above.

\begin{proof}[Proof of Theorem \ref{thm:rateofconvPolynom}]
As above, we choose $p>(1-\epsilon)/\epsilon$ large enough and apply Theorem \ref{thm:SmoothingThm} to $\mu=\mu_n^W,\nu=\mu_\infty, K=2,\eta=1$, $z^*=0$ and $a=n^{1/2-\epsilon}$, and obtain
 \begin{align*}
 D(\mu_n^W,\mu_\infty)\lesssim &\ n^{1/2}\norm {U_n^W-U_\infty+1/2}_{L^p(B_{3}(0)\setminus B_{n^{-1/2+\epsilon}}(0))}\\
 &+\mu_n^W(B_{2n^{-1/2+\epsilon}}(0))+\mu_n^W(B_{2}(0)^c)\\
 &+\mu_\infty(B_{2n^{-1/2+\epsilon}}(0))+\mu_\infty(B_{2}(0)^c)\\
 &+\sup_{R\geq 0,z_0\in\IC}\mu_\infty\left( R\leq \abs{\cdot-z_0 }\leq R+2n^{-1/2+\epsilon}\right).
 \end{align*}
Let us consider each term starting with the last one. Obviously the last term is of order $n^{-1/2+\epsilon}$ and the third line equals $4\pi n^{-1+2\epsilon}$. From an already existing (non-uniform) Local Circular Law for random polynomials, see \cite{TV15Polynomials} formula (87), it follows that with overwhelming probability the second line of our estimation can also be bounded by $cn^{-1+2\epsilon}$. Therefore it remains to control the $L^p$ distance of the logarithmic potentials. The application of Monte Carlo sampling and the pointwise control of the logarithmic potentials from \eqref{eq:ConcLogPotPoly} remains unchanged. The only notable difference to the proof of Theorem \ref{thm:rateofconv} is the restriction to polynomially bounded moduli of the zeros. From Rouch\'e's Theorem, we deduce an upper bound for the largest root
\begin{align*}
\abs{\lambda}_{\max}\le 1+\frac{\max\{c_0\abs{\xi_0},\dots,c_{n-1}\abs{\xi_{n-1}}  \}}{c_n\abs{\xi_n}}
\end{align*}
of any polynomial. Hence for any $Q>0$ we have
\begin{align*}
\IP(\abs\lambda_{\max}\ge n^{(Q+1)/\delta})&\le \IP\left(\frac{\max\{\abs{\xi_0},\dots,\abs{\xi_{n-1}} \}}{\abs{\xi_n}}\gtrsim n^{(Q+1)/\delta}\right)\\
&\le (n-1)\IP(\abs{\xi_0}\gtrsim n^{(Q+1)/\delta}\abs{\xi_n}) \\
&\lesssim \frac{n-1}{n^{Q+1}}\IE\abs{\xi_0}^\delta\IE\abs{1/\xi_0}^\delta\lesssim n^{-Q},
\end{align*}
which replaces \eqref{eq:s_max} and the proof is finished.
\end{proof}

Lastly, we provide the elementary 

\begin{proof}[Proof of Lemma \ref{lem:GinibreRate}]
Since \cite{Gin65}, the density $p_n$ of $\bar\mu_n$ has been known to be
\begin{align}\label{eq:DensityGinibre}
p_n(z)=\frac 1 \pi e^{-n\abs z ^2} \sum_{k=0}^{n-1}\frac{n^k\abs z ^{2k}}{k!},
\end{align}
which converges to $p_\infty(z)=\frac 1 \pi\eins_{B_1(0)}(z)$. In the case of $z_0=0$, we can explicitly calculate
\begin{align*} 
 \bar\mu_n(B_R(0))&=\frac 1 \pi\int_{B_R(0)}e^{-n\abs z ^2}\sum_{k=0}^{n-1}\frac{n^k\abs z ^{2k}}{k!}dz\\
 &= \frac 1 n\sum_{k=0}^{n-1}  \int_0^{nR^2}e^{-r}\frac {r^k}{k!} dr\\
 &=\frac 1 n\sum_{k=0}^{n-1} 1-e^{-nR^2}\sum_{j=0}^k \frac{(nR^2)^j}{j!}\\
 &=1-e^{-nR^2}\sum_{k=0}^{n-1}  \frac{(n-k)(nR^2)^k}{nk!}\\
 &=1-e^{-nR^2}\left( \frac{(nR^2)^{n}}{n!}+(1-R^2)\sum_{k=0}^{n-1}  \frac{(nR^2)^k}{k!}\right),
\end{align*}
where we used the substitution $r=n\abs z ^2$ and integration by parts. The function
\begin{align*}
\bar D_n(R)&=\mu_\infty (B_R(0)) -\bar\mu_n(B_R(0))\\
&=1\wedge R^2-1+e^{-nR^2}\left( \frac{(nR^2)^{n}}{n!}+(1-R^2)\sum_{k=0}^{n-1}  \frac{(nR^2)^k}{k!}\right)
\end{align*}
is continuous in $R$ and differentiable for $R\neq 1$ with radial derivative as above
\begin{align*}
 2R\left(\eins_{[0,1)}(R) -e^{-nR^2} \sum_{k=0}^{n-1}\frac{(nR^2)^k}{k!} \right)\begin{cases}                                                        >0\quad&\text{, if }R<1\\                                                                               <0\quad&\text{, if }R>1
\end{cases}.
\end{align*}
Hence the maximum is attained at $R=1$ and Stirling's formula yields
\begin{align}\label{eq:AsymptoticGinibreStirling}
 \sup_{R>0}\abs{\bar\mu_n(B_R(0))-\mu_\infty (B_R(0)) }=\mu_\infty (B_1(0))-\bar\mu_n (B_1(0))=\frac{n^{n}}{e^nn!}\sim \frac 1{\sqrt{2\pi n }}.
\end{align}
The distances of arbitrary balls are likewise bounded by
\begin{align*}
\mu_\infty (B_R(z_0))-\bar\mu_n (B_R(z_0))&\leq \mu_\infty (B_1(0))-\bar\mu_n (B_1(0))\\
\bar\mu_n (B_R(z_0))-\mu_\infty (B_R(z_0))&\leq \bar\mu_n (B_1(0)^c)=\mu_\infty (B_1(0))-\bar\mu_n (B_1(0)),
\end{align*}
hence the first part of the statement is proven. For $R\leq1$ we have
\begin{align*}
 \bar D_n(R)=e^{-nR^2}\left( \frac{(nR^2)^n }{n!}-(1-R^2)\sum_{k=n}^\infty \frac{(nR^2)^k }{k!}\right)
\end{align*}
and
\begin{align*}
e^{-nR ^2} \sum_{k=n}^{\infty}\frac{n^kR ^{2k}}{k!}&\leq  e^{-nR ^2} \frac{(nR^2)^n}{ n!}\sum_{k=0}^\infty\left( \frac{nR ^{2}}{(n+1)}\right)^k\\
&=e^{-nR ^2} \frac{(nR^2)^n}{n!} \frac{n+1}{n(1-R ^2)+1}\\
&\sim\frac 1 {\sqrt{2\pi n}} e^{-n(R^2 -1-\log (R ^2))}\frac{n+1}{n(1-R ^2)+1},
\end{align*}
where we applied Stirling's formula again. Consequently
\begin{align*}
\abs{\bar D_n(R)}&\lesssim\frac{1}{\sqrt{n}} e^{-n(R^2 -1-\log (R ^2))}\left(1+(1-R^2)\frac{n+1}{n(1-R ^2)+1}\right)\\
&\lesssim\frac{1}{\sqrt{n}} e^{-n(R^2 -1-\log (R ^2))}
\end{align*}
for $R\leq 1$. On the other hand if $R\geq 1$, then
\begin{align*}
 \bar D_n(R)=e^{-nR^2}\left( \frac{(nR^2)^n }{n!}-(R^2-1)\sum_{k=0}^{n-1} \frac{(nR^2)^k }{k!}\right),
\end{align*}
where analogously we have
\begin{align*}
\sum_{k=0}^{n-1} \frac{(nR^2)^k }{k!}\le\frac{(nR^2)^{n-1}}{(n-1)!}\sum_{k=0}^{n-1}\left(\frac{n-1}{nR ^2} \right)^k\le \frac{(nR^2)^{n}}{(n)!} \frac{1}{(R ^2-1)+1}
\end{align*}
and hence
\begin{align*}
\abs{\bar D_n(R)}&\lesssim\frac{1}{\sqrt{n}} e^{-n(R^2 -1-\log (R ^2))}.
\end{align*}
Finally choose $R=1-\epsilon$ (or $R=1+\epsilon$, respectively) and note that $R^2-1-\log R^2\geq 2\epsilon^2+\mathcal O(\epsilon ^3)$, we conclude
\begin{align*}
\abs{\bar D_n(1-\epsilon)}\lesssim e^{-n\epsilon^2}
\end{align*}
and the second part of the Lemma follows.
\end{proof}

\appendix
\renewcommand{\thesection}{\Alph{section}}
\section{Appendix}\label{sec:Appendix}

\begin{lemma}\label{lem:Convergence}
Convergence of distributions on $\IC$ with respect to the spherical Kolmogorov distance $ D$ implies weak convergence.
\end{lemma}
For absolutely continuous limit distributions, the converse statement is also true, see for instance \cite{Dud76}. Hence $ D$ is a reasonable object for studying the rate of convergence to the Circular Law. Moreover we justify the term \emph{Kolmogorov distance} by formally retrieving the 1-di\-men\-sio\-nal Kolmogorov distance $ d(\mu_j,\nu_j)$ of the marginals $j=1,2$ in limits such as
\begin{align*}
(\mu_1-\nu_1)((-\infty,t])=\lim_{K\to \infty} (\nu-\mu)(B_K(t+K,0)).
\end{align*}

\begin{proof}
We will prove vague convergence and tightness. Let $\mu,\nu$ be distributions on $\IC$, $f\in\mathcal C _c(\IC)$ be a continuous function with compact support and $f_r=\frac 1{\pi r^2}f\ast\eins_{B_r(0)}$ be its ball mean function. Furthermore denote by $\lebesgue$ the Lebesgue measure on $\IC$ and set $\eta=\mu-\nu+\lebesgue$. It holds
\begin{align*}
\int f d(\mu-\nu)-\int (f-f_r)d\eta&=\int f_rd\eta-\int fd\lebesgue\\
&=\int \int f(y)\frac 1{\pi r ^2}\eins_{B_r(0)} (y-x)d\lebesgue(y)d\eta(x)-\int fd\lebesgue\\
&=\int f(y)\left(\int\frac 1{\pi r ^2}\eins_{B_r(y)} (x)d\eta(x)    -1\right) d\lebesgue(y)\\
&=\frac1{\pi r ^2}\int f(y)\left(\mu(B_r(y))-\nu(B_r(y))\right) d\lebesgue(y).
\end{align*}
Now choosing a sequence $\mu=\mu_n$ converging to $\nu$ with respect to $ D$ implies for all $r>0$
\begin{align*}
\abs{\int f d(\mu_n-\nu)}&\le \abs{\int (f-f_r)d\eta_n}+\frac1{\pi r ^2}\int f(y)\abs{\mu(B_r(y))-\nu(B_r(y))} d\lebesgue(y)\\
&\le \int \abs{f-f_r}d(\mu_n+\nu+\lebesgue)+\frac{ D(\mu_n,\nu)}{\pi r ^2}\int \abs{f(y)} d\lebesgue(y)\\
&\le 2\norm{f-f_r}_{L^\infty(\lebesgue)} +\norm{f-f_r}_{L^1(\lebesgue)}+\frac{ D(\mu_n,\nu)}{\pi r ^2}\int \abs{f(y)} d\lebesgue(y).
\end{align*}
First as $n\to\infty$, the last term converges to $0$, then as $r\to 0$, the first term vanishes due to the continuity of $f$ and the second due to Lebesgue Differentiation Theorem. This implies that $\mu_n$ converges to $\nu$ in weak$\ast$ convergence.

Tightness follows easily from the convergence in $D$. For any $\epsilon$ let $N\in\IN$ be sufficiently large for $D(\mu_n,\mu)\le\epsilon/2$ for all $n>N$. Then choose $K_\epsilon>0$ sufficiently large for $\mu(B^c_{K_\epsilon}(0)),\mu_n(B^c_{K_\epsilon}(0))\le\epsilon/2$ for all $n\le N$, then $\mu_n(B^c_{K_\epsilon}(0))\le \epsilon$ for all $n$.
\end{proof}

In order to prove Proposition \ref{thm:ConcLogPot}, we will directly follow the approach of \cite{GNT17Local}, making use of Girko's Hermitization trick to convert the non-Hermitian problem into a Hermitian one, apply the local Stieltjes transform estimate from \cite{GNT17Local} and the smoothing inequality from \cite{GT03RateSemicircular}. Let $\tilde\nu^z_n$ be the symmetrized empirical singular value distribution of the shifted matrices $X/\sqrt n -z$, defined in \eqref{eq:GirkoHerm} and
\begin{align*}
 m_n(z,\cdot ):\IC\setminus\IR\to\IC, w\mapsto\int_{\IR}\frac{1}{w-t}d\tilde\nu^z_n(t)
\end{align*}
be the Stieltjes transform which converges a.s. to the solution of 
\begin{align}\label{eq:StieltjesEquation}
s(z,w)=-\frac{s(z,w)+w}{(w+s(z,w))^2-\abs z ^2}
\end{align}
see for instance \cite{GT10Circular}. It is known that $s(z,\cdot )$ corresponds to a limiting measure $\tilde \nu ^z$ which has a symmetric bounded density $\rho^z$ (the bound holds uniformly in $z$) and has compact support
\begin{align*}
 \mathbb{J}^z:=\begin{cases}
              [-\lambda_+,-\lambda_-]\cup[\lambda_-,\lambda_+], &\text{ if }\abs z>1\\
              [-\lambda_+,\lambda_+], &\text{ if }\abs z\leq 1
             \end{cases},
\end{align*}
where the endpoints are given by
\begin{align*}
 \lambda_\pm^2:=\frac{(\alpha\pm 3)^3}{8(\alpha\pm 1)}\wedge 0,\quad \alpha:=\sqrt{1+8\abs z^2}.
\end{align*}
Note that $\lambda_-\sim (1-\abs z )^{3/2}$ as $\abs z\to 1$, i.e. a new gap in the support emerges at $0$. Therefore $s$ will be unbounded for $z$ close to the edge, which is the reason for the bulk constraint of Proposition \ref{thm:ConcLogPot}. 

\begin{proof}[Proof of Proposition \ref{thm:ConcLogPot}]
Fix some arbitrary $Q,\tau>0$ and $z\in B_{1+\tau^{-1}}(0)$ satisfying $\abs{1-\abs z }\ge \tau$. As is explained in Girko's Hermitization trick \eqref{eq:GirkoHerm},
\begin{align*}
 \abs{U_n(z)-U_\infty(z) }=\abs{\int_{\IR} \log \abs x d(\tilde\nu^z_n-\tilde\nu^z)(x)	}
\end{align*}
and therefore it is necessary to estimate the extremal singular values as well as the rate of convergence of $\tilde\nu^z_n$ to $\tilde\nu^z$ in Kolmogorov distance $ d_n^*(z)$. Introduce the events
\begin{align*}
 \Omega_0:=\{s_{\min }\geq n^{-B}  \},\quad \Omega_1:=\{s_{\max }\leq n^{B'}  \},\quad \Omega_2:=\{  d_n^*(z)\leq c \log^3n/n\}
\end{align*}
for some constants $B,B',c>0$ yet to be chosen. Theorem 2.1 in \cite{TV08Circular} states that there exists a constant $B>0$ such that $\IP(\Omega_0^c)\lesssim n^{-Q}$ and analogously to what has been shown in \eqref{eq:s_max} there exists a constants $B'>0$ with $\IP(\Omega_1^c)\lesssim n^{-Q}$. Since $\tilde\nu^z$ has a bounded density, we get
\begin{align*}
 \abs{\int_{-n^{-B}}^{n^{-B}}\log \abs x d\tilde\nu^z(x)}\lesssim \log n n^{-B}
\end{align*}
and furthermore on $\Omega_2$ it holds that
\begin{align*}
 \abs{\int_{n^{-B}\leq \abs x\leq n^{B'}}\log\abs xd(\tilde\nu_n^z-\tilde\nu ^z )(x)}\lesssim   d^*_n(z)\log n\lesssim \frac{\log ^4n}{n}.
\end{align*}
Hence the claimed concentration of $U_n$ holds on $\Omega_0\cap\Omega_1\cap\Omega_2$, implying
\begin{align*}
 \IP\left( \abs{U_n(z)-U_\infty(z) } \geq c \frac{\log ^4 n}{n}\right)\leq \IP(\Omega_0^c)+\IP(\Omega_1^c)+\IP(\Omega_2^c)
\end{align*}
and it remains to check $\IP(\Omega_2^c)\leq n^{-Q}$, which has been done explicitly in \cite{GNT17Local}, (4.14)-(4.16) using the smoothing inequality [Corollary B.3] from \cite{GT03RateSemicircular} and the local law for $ d_n^*(z)$ in terms of their Stieltjes transforms.
\end{proof}

Very recently, in \cite{AEK19}, the restriction to the bulk has been removed by the usage of an cusp fluctuation averaging method. For the proof of Proposition \ref{prop:ConcLogPotAEK}, we will use the following identity that goes back to Tao and Vu \cite{TV15uni}.

\begin{lemma}\label{lem:logdetRepr} For any $T>0$ it holds
\begin{align}\label{eq:logdetRepr}
 U_n(z)=\int_0^{T} \Im(m_n(z,i\eta))d\eta-\frac 1 {2n}\log\abs{\det (V(z)-iT)}
 \end{align}
and in the limit $T\to\infty$ we also have
\begin{align}\label{eq:logpotRepr}
 U_\infty(z)=\int_0^{\infty}\Im (s(z,i\eta)) -\frac{1}{1+\eta}d\eta.
\end{align}

 \end{lemma}
 \begin{proof}
The distributional equation for Stieltjes transforms $m_\mu=2\partial_z U_\mu$ in terms of Wirtinger derivatives the case of $\mu=\tilde\nu_n^z$ yields $\Im(m_n(z,i\eta))=-\partial_\eta U_{\tilde\nu_n^z}(i\eta)$. We integrate with respect to $\eta$ and readily obtain
\begin{align*}
 \int_0^{T}\Im (m_n(z,i\eta))d\eta=-U_{\tilde\nu_n^z}(iT)+U_{\tilde\nu_n^z}(0).
\end{align*}
The first term is as in the claim and the second term follows from rephrasing Girko's Hermitization Trick \eqref{eq:GirkoHerm} as
 \begin{align*}
  U_n(z)=U_{\mu_n}(z)=U_{\tilde\nu_n^z}(0)=U_{\nu_n^z}(0).
 \end{align*}
The same arguments yield 
\begin{align*}
 U_\infty(z)&=\int_0^T \Im (s(z,i\eta))d\eta -\int_{-\lambda_+}^{\lambda_+}\log \abs{x-iT}d\tilde\nu_\infty^z(x)\\
 &=\int_0^T \Im (s(z,i\eta))d\eta-\int_0^T\frac{1}{1+\eta}d\eta+\mathcal O(T^{-1}),
\end{align*}
since the second integral is asymptotically equivalent to $\log T$ as $T\to \infty$. For large $\eta$ we see that
\begin{align}\label{eq:etaSquare}
 \Im (s(z,i\eta))-\frac 1 {1+\eta}=\int\frac{1-x^2/\eta}{\left(x^2/\eta+\eta\right)\left(1+\eta\right)}d\tilde\nu^z_\infty(x)\sim \frac 1 {\eta^2},
\end{align}
hence it is integrable and we may pass to the limit $T\to \infty$ to obtain the claim.
 \end{proof}

 \begin{proof}[Proof of Proposition \ref{prop:ConcLogPotAEK}]
We will show all steps until a result from \cite{AEK19} can be directly applied. Using Lemma \ref{lem:logdetRepr}, we need to estimate
\begin{align}\label{eq:(6.1)}
U_\infty(z)-U_n(z)&=\int_0^{T}\Im (s(z,i\eta) - m_n(z,i\eta))d\eta\nonumber\\
&+\int_T^{\infty}\Im (s(z,i\eta)) -\frac{1}{1+\eta}d\eta\nonumber\\
&+\frac 1 {2n}\log\abs{\det (V(z)-iT)}-\int_0^{T}\frac{1}{1+\eta}d\eta,
\end{align}
which corresponds to a pointwise estimate of the integral in \cite[Equation (6.1)]{AEK19}. By \eqref{eq:etaSquare}, the second term is of order $\mathcal O (T^{-1})$. Regarding the third term it holds
\begin{align}\label{eq:V-inc}
 \frac 1 {2n}\log\abs{\det (V(z)-iT)}&=\log T    +\log\abs{\det( i-T^{-1}V(z))}\nonumber\\&=\log T + \sum_{j=1}^n\log\left( 1+\frac{s_j(z)^2}{T^2}\right),
\end{align}
where $s_j(z)$ are the non-negative eigenvalues of $V(z)$, or equivalently the singular values of $X/\sqrt n-z$. Similar to \eqref{eq:s_max}, we have $
 \IP(s_{\max }(z)\geq n^{(Q+1)/2})\lesssim n^{-Q}$. Thus, the last term in \eqref{eq:V-inc} is neglegible if we choose $T=n^C$ for $C$ large enough. The remaining $\log T$ cancels the last term of \eqref{eq:(6.1)}. Again, by \cite[Theorem 2.1]{TV08Circular} it holds $\IP(s_{\min }(z)< n^{-B})\lesssim n^{-Q}$ and consequently, Equation \eqref{eq:(6.1)} becomes 
\begin{align*}
 U_\infty(z)-U_n(z)&=\int_{n^{-B}}^{n^C}\Im (s(z,i\eta)-m_n(z,i\eta))d\eta +\mathcal O (n^{-1})
\end{align*}
with overwhelming probability. At this stage, according to \cite[Remark 6.2]{AEK19}, it holds \cite[Lemma 6.1]{AEK19} stating 
\begin{align*}
 \IE\Big|\int_{n^{-B}}^{n^C}\Im (s(z,i\eta)-m_n(z,i\eta))\Big|^p\lesssim \frac {n^{\delta p}}{n^p}
\end{align*}
for any $\delta>0$, $p\in\IN$. Choosing $\delta=\epsilon/2$ and $p$ sufficiently large, an application of Markov's inequality finishes the proof. We shall point out that the conditions \cite[(6.4) Remark 2.5]{AEK19} are satisfied in our case, while \cite[(A1), (A2)]{AEK19} coincide with condition (A) in our claim.
 \end{proof}

\section*{Acknowledgements}
Financial support by the German Research Foundation (DFG) through the IRTG 2235 is gratefully acknowledged. We would like to thank A. Tikhomirov and A. Naumov for helpful discussions and valuable suggestions. Furthermore, we thank the referee for the useful feedback and remarks.

\end{document}